\documentclass[11pt,a4paper,reqno]{amsart}
\usepackage[T1]{fontenc}
\usepackage{lmodern}
\usepackage[a4paper,margin=28mm]{geometry}
\usepackage{amsmath,amssymb}
\usepackage{xcolor}
\definecolor{LinkBlue}{HTML}{002FA7}
\definecolor{CitationRed}{HTML}{800020}
\definecolor{UrlGray}{HTML}{304C68}
\usepackage[colorlinks=true,linkcolor=LinkBlue,citecolor=CitationRed,
  urlcolor=UrlGray,bookmarksnumbered=true]{hyperref}
\hypersetup{pdftitle={Good minimal models and termination with scaling for compact Kähler pairs}}
\newtheorem{theorem}{Theorem}[section]
\newtheorem{proposition}[theorem]{Proposition}
\newtheorem{lemma}[theorem]{Lemma}
\newtheorem{corollary}[theorem]{Corollary}
\theoremstyle{definition}
\newtheorem{gzdstatement}{}

\newtheorem{example}[theorem]{Example}
\theoremstyle{plain}
\numberwithin{equation}{section}
\newcommand{\R}{\mathbb R}
\newcommand{\Q}{\mathbb Q}

\newcommand{\Z}{\mathbb Z}
\newcommand{\M}{\mathbf M}
\newcommand{\cO}{\mathcal O}
\DeclareMathOperator{\Exc}{Exc}
\DeclareMathOperator{\Supp}{Supp}
\DeclareMathOperator{\coeff}{coeff}
\title[Good minimal models of compact K\"ahler pairs]
  {Good minimal models and termination with scaling for compact K\"ahler pairs}
\author{Zhengyu Hu}
\address{Mathematical Sciences Research Center, Chongqing University of Technology,
Chongqing 400054, China}
\email{zhengyuhu16@gmail.com}
\author{Kefeng Liu}
\address{Mathematical Sciences Research Center, Chongqing University of Technology,
Chongqing 400054, China}
\email{kefengliu@cqut.edu.cn}
\hypersetup{pdfauthor={Zhengyu Hu and Kefeng Liu}}
\date{November 30, 2026}
\begin{document}

\begin{abstract}
The good divisorial Zariski decomposition (GZD) asserts that if
$(X,B)$ is a log smooth compact K\"ahler pair with boundary $\Q$-divisor
$B$ and $K_X+B$ is pseudoeffective, then on some modification
$\mu\colon W\to X$ the divisorial Zariski decomposition
$\mu^*(K_X+B)\sim_\Q P+N$ has semiample positive part $P$.
This was proved by OpenAI in a preprint. Assuming (GZD), we prove that
a log canonical pair $(X,B)$ on a normal compact K\"ahler space has a
good minimal model or a Mori fibre space. We also prove that
any minimal model program for $K_X+B$ with K\"ahler scaling terminates
if $(X,B)$ is a divisorial log terminal pair and $X$ is globally strongly
$\Q$-factorial. In both results the boundary $B$ may have real coefficients.
\end{abstract}

\subjclass[2020]{Primary 14E30; Secondary 32J27, 32J25}
\keywords{Minimal model program, compact K\"ahler spaces, good minimal models,
termination with scaling, divisorial Zariski decomposition}

\maketitle

\section{Introduction}

We assume \ref{gzd} throughout. In \cite[Theorem~2.13]{KA}, OpenAI
proved the statement using projective logarithmic Iitaka subadditivity
(\cite[Theorem~1.2]{PA}). For a pseudoeffective real $(1,1)$-class
$\alpha$ on a compact complex manifold, the divisorial negative part is
$N_\sigma(\alpha)=\sum_D\nu_D(\alpha)D$, where $D$ runs over prime divisors
and $\nu_D(\alpha)$ is the minimal multiplicity along $D$
\cite[\S3]{Bou}. We also use snc for simple normal crossings.

\begin{gzdstatement}\label{gzd}
Let $X$ be a smooth compact K\"ahler manifold and $B$ a boundary
$\Q$-divisor on $X$ with snc support. If $K_X+B$ is pseudoeffective,
then there is a modification $\mu\colon W\to X$ from a smooth compact
K\"ahler manifold $W$ such that
$$
 \mu^*(K_X+B)\sim_\Q P+N,\qquad
 N=N_\sigma\big(\mu^*(K_X+B)\big),
$$
where $P$ is a semiample $\Q$-line bundle and $N$ is an effective
$\Q$-divisor.
\end{gzdstatement}

Under \ref{gzd} we prove two theorems for compact K\"ahler pairs whose
boundaries may have real coefficients. For log canonical (lc) pairs we
obtain good minimal models or Mori fibre spaces
(Theorem~\ref{main}). For divisorial log terminal (dlt)
pairs on globally strongly $\Q$-factorial spaces we prove termination of
any minimal model program (MMP) for $K_X+B$ with K\"ahler scaling
(Theorem~\ref{scaling}).

\begin{theorem}\label{main}
Let $(X,B)$ be an lc pair on a normal compact K\"ahler space.
Then $(X,B)$ has either a good minimal model $(Y,B_Y)$ or a Mori fibre
space $(Y,B_Y)\to Z$. In both cases $Y$ is compact K\"ahler. Moreover,
if $(X,B)$ is dlt, then we can take $X\dashrightarrow Y$ to be a
bimeromorphic contraction.
\end{theorem}

\begin{corollary}\label{nef}
Let $(X,B)$ be an lc pair on a normal compact K\"ahler space.
Suppose that $K_X+B$ is nef. Then $K_X+B\sim_\R g^*H$ for a
contraction $g\colon X\to Z$ onto a normal projective variety $Z$ and an
ample $\R$-divisor $H$ on $Z$. Moreover, if $B$ is a $\Q$-divisor,
then $m(K_X+B)$ is Cartier and base point free for some positive
integer $m$.
\end{corollary}

In particular, Corollary~\ref{nef} proves part~(1) of the abundance conjecture
of C.~D.~Hacon and L.~Xie \cite[Conjecture~1.16]{HX}.

\begin{theorem}\label{scaling}
Let $(X,B)$ be a dlt pair on a globally strongly $\Q$-factorial
compact K\"ahler space. Let $\omega$ be a K\"ahler class such that
$c_1(K_X+B)+\omega$ is nef. Then we can run a $(K_X+B)$-MMP with
scaling of $\omega$, and any such MMP terminates. If $K_X+B$ is
pseudoeffective, then the final model is a good minimal model of
$(X,B)$. Otherwise the MMP ends with a Mori fibre space.
\end{theorem}

For a reduced snc divisor $\sum_i S_i$ on a smooth compact K\"ahler
manifold $T$, the spaces $H^0(T,dK_T+\sum_i\alpha_i S_i)$,
$0\le\alpha_i\le d$, are the weight spaces of one log canonical ring.
By \ref{gzd}, a Veronese subring is finitely generated. From its
weights we obtain a rational polyhedral decomposition and extend the
divisorial Zariski decomposition to real coefficients on each cell (Proposition~\ref{realdecomp}).

The rest of the proof follows the earlier approach of the first author.
In \cite[proof of Lemma~2.8]{Hu20} the first author compared the
coefficients of a vertical divisor after equidimensional reduction. Here
we make this comparison using Gauduchon degrees on the fibre components
(Lemma~\ref{very-exceptional}). We remark that it is not enough to test on
fibre curves (Example~\ref{fibre-curve-test}).
As in \cite[Lemma~3.14]{Hu25}, we find a model on which the transform of
$K_X+B+s\omega$ is nef for sufficiently small $s\ge0$. We compare with
the current MMP model at two distinct consecutive scaling numbers in
this range and conclude that it is minimal (Lemma~\ref{interval}).

\noindent\textbf{Related work.}
Under a good decomposition hypothesis, Y.~Li and S.~Rao construct a
finite $K_X$-MMP when $X$ is a compact K\"ahler manifold and $K_X$ is
pseudoeffective \cite[Theorem~1.4]{LR}. C.~D.~Hacon and P.~Wang consider
Kawamata log terminal (klt) pairs $(X,B)$ on compact K\"ahler spaces,
where $B$ is a $\Q$-divisor. When $X$ is strongly $\Q$-factorial,
they prove that any $(K_X+B)$-MMP with K\"ahler scaling terminates
\cite[Theorem~1.1]{HW}. In the pseudoeffective case they use OpenAI's
theorem on divisorial valuations for finite cyclic group actions to bound
the Cartier indices along the MMP \cite[Theorem~4.1 and Proposition~5.2]{HW}.
On the minimal model they apply OpenAI's abundance theorem
\cite[Theorem~1.2]{KA} to obtain a good minimal model
\cite[Corollary~1.2]{HW}.
We obtained our results independently of Li--Rao and Hacon--Wang.
For projective pairs, see \cite[Theorem~1.1]{BH} for minimal models
from Zariski decompositions and \cite[Theorem~4.1(iii)]{Bir} for
termination with scaling. In the K\"ahler setting, O.~Das,
C.~D.~Hacon and M.~P\u{a}un treated $\Q$-factorial dlt pairs of dimension four with
$K_X+B\sim_\Q E\geq0$ \cite[Theorem~1.1]{DHP}.

\noindent\textbf{AI disclosure.}
The mathematical ideas, the proofs and all details of the proofs are due to the authors. Under the authors' direction, AI tools were used to draft the text. The authors revised the manuscript, verified all results, and take full responsibility for it.

\section{Preliminaries}

\noindent\textbf{Pairs and classes.}
All complex spaces are irreducible, and for pairs, nef classes and
pseudoeffective classes we use the definitions of \cite[\S2]{DHP}. In
particular, for a pair $(X,B)$ we require $X$ to be normal, $B$ a boundary
$\R$-divisor and $K_X+B$ $\R$-Cartier. Canonical divisors on different
models are compared by local canonical forms, and log resolutions can be
taken projective.

The real Bott--Chern group $H^{1,1}_{\mathrm{BC}}(X,\R)$ consists of
closed real $(1,1)$-forms with smooth local potentials, modulo
$i\partial\bar\partial u$ for global smooth real functions $u$. A class
is pseudoeffective if its pullback to a resolution contains a closed
positive $(1,1)$-current. We write $L$ for $c_1(L)$ in equalities of
classes.

Recall that a class $\alpha$ on a reduced compact complex space $F$ with a
Hermitian form $\gamma$ is nef if for each $\varepsilon>0$ it has a smooth
representative $\theta_\varepsilon\geq-\varepsilon\gamma$. By nefness on the
fibres we mean nefness on each reduced fibre.

\medskip\noindent\textbf{Semiample line bundles.}
An $\R$-line bundle $L$ is semiample if $L=\sum_i a_iL_i$ in
$\operatorname{Pic}(X)\otimes_\Z\R$, where $a_i>0$ and each $L_i$ is a
globally generated line bundle.

\begin{lemma}\label{semiample}
Let $L$ be an $\R$-line bundle on a normal compact complex space $X$.
Then $L$ is semiample if and only if $L=g^*A$ for a contraction
$g\colon X\to Z$ onto a normal projective variety $Z$ and an ample
$\R$-line bundle $A$ on $Z$. Moreover, if $L$ is a semiample
$\Q$-line bundle, then $mL$ is a base point free line bundle for
sufficiently divisible $m>0$.
\end{lemma}
\begin{proof}
Suppose $L=\sum_i a_iL_i$. Let $h\circ g$ be the Stein factorization
of the product of the maps given by the $|L_i|$. Since $h$ is finite,
$A=h^*(\sum_i a_i\,\mathrm{pr}_i^*\cO(1))$ is ample. We have $L=g^*A$.
Conversely, an ample $\R$-line bundle is a positive real
combination of very ample ones. Since $g^*$ is injective on Picard groups
(projection formula), for a $\Q$-line bundle $L$ we can take $A$ ample
$\Q$-Cartier. Then $mL$ is base point free for sufficiently divisible $m$.
\end{proof}

\medskip\noindent\textbf{Negativity lemma.}
Let $f\colon X\to Y$ be a proper surjection with connected fibres. A
divisor $D$ on $X$ is \emph{vertical} if $f(\Supp D)\subsetneq Y$, and
it is \emph{very exceptional} if moreover every prime divisor on $Y$ is
the image of a prime divisor on $X$ not contained in $\Supp D$
\cite[Definition~3.1]{Bir}.

\begin{lemma}[Negativity for vertical divisors]\label{very-exceptional}
Let $f\colon X\to Y$ be a proper surjection with connected fibres
between normal compact complex spaces. Let $D=D^+-D^-$ be a vertical
$\R$-Cartier divisor on $X$, where $D^+,D^-\geq0$ have no common
components. Suppose that $c_1(D)$ is nef on the fibres of $f$. Then
the following hold.
\begin{enumerate}
\item If $D^+$ is very exceptional over $Y$, then $D\leq0$.
\item Let
$$
 \begin{array}{ccc}
 X' & \xrightarrow{\ p\ } & X \\
 {\scriptstyle f'}\downarrow && \downarrow{\scriptstyle f}\\
 Y' & \xrightarrow{\ q\ } & Y
 \end{array}
$$
be an equidimensional reduction of $f$, where $Y'$ is smooth and $p,q$
are projective modifications. Then there is a $q$-nef $\R$-divisor
$D_{Y'}$ on $Y'$ such that $p^*D=f'^*D_{Y'}$.
\end{enumerate}
\end{lemma}
\begin{proof}
\emph{Step 1. Descent to $Y'$.} We obtain the diagram in (2) by
analytic flattening and a resolution of the base
\cite[proof of Lemma~2.1]{DO}. The fibres of $f'$ are connected (Stein
factorization). Since $f'$ is equidimensional, each component of
$p^*D$ dominates a prime divisor on $Y'$. As in
\cite[proof of Lemma~2.8]{Hu20}, for each such prime divisor $Q$ we write
$$
 f'^*Q=\sum_i m_iR_i,\qquad d_i=\coeff_{R_i}(p^*D),\qquad
 d_Q=\max_i d_i/m_i .
$$
Near the fibre over a general point $y\in Q$, the divisor
$G=d_Qf'^*Q-p^*D$ is effective and omits at least one $R_i$.
Since $R_i\cap R_j$ has smaller fibre dimension over $Q$ for $i\ne j$,
each component of the fibre is contained in exactly one $R_i$.
If $G$ contains some $R_j$, by connectedness a component
$F\not\subset\Supp G$ intersects $\Supp G$. Then $G|_U$ is effective
and nonzero on a resolution $U$ of $F$.

Let $d=\dim U$ and $\gamma$ be a Gauduchon metric on $U$ \cite{Gau}.
We integrate against $\gamma^{d-1}$ to define a degree on Bott--Chern
classes ($\partial\bar\partial\gamma^{d-1}=0$). The class of $f'^*Q$
vanishes on $U$. Since $U$ maps into a fibre of $f$, $c_1(p^*D)|_U$ is
nef. Its degree is nonnegative (integrate representatives
$\geq-\varepsilon\gamma$ and let $\varepsilon\to0$). We get
$$
 0<\int_{G|_U}\gamma^{d-1}
   =-\deg_\gamma(c_1(p^*D)|_U)\leq0 ,
$$
a contradiction. Hence $d_i=m_id_Q$ for each $i$, or equivalently
$p^*D=f'^*D_{Y'}$ with $D_{Y'}=\sum_Q d_Q Q$.

\emph{Step 2. Nefness over $Y$.} For a curve $C$ contracted by $q$,
let $U$ resolve a component of $f'^{-1}(C)$ dominating $C$ and let
$h\colon U\to C^\nu$ be the induced map. A $(1,1)$-class on $C^\nu$
is determined by its degree. For a positive form $\omega$ of integral
one on $C^\nu$ and a Gauduchon metric $\gamma$ on $U$ we have
$$
 \deg_\gamma(c_1(p^*D)|_U)
  =(D_{Y'}\cdot C)\int_U h^*\omega\wedge\gamma^{\dim U-1}.
$$
The degree is nonnegative by nefness on the fibres. Since $h$ is
surjective, the integral is positive and $D_{Y'}\cdot C\ge0$.
By projectivity of $q$, $D_{Y'}$ is $q$-nef.

For (1), let $Q$ be a prime divisor not contracted by $q$. Since
$D^+$ is very exceptional, a prime divisor over $q(Q)$ has coefficient
$\le0$ in $D$. On its strict transform $R_i$ we get $d_Q=d_i/m_i\le0$.
Thus $q_*D_{Y'}\le0$. By the negativity lemma
\cite[Lemma~1.3]{Wang}, \cite[Lemma~2.5]{HX} we have $D_{Y'}\le0$ and
$D=p_*f'^*D_{Y'}\le0$.
\end{proof}

\begin{example}\label{fibre-curve-test}
Both conclusions of Lemma~\ref{very-exceptional} can fail if nefness is
tested only on fibre curves. Let $\mathbb T$ be a
two-dimensional torus without curves, let $X$ be the blowup of
$\mathbb T\times\mathbb P^1$ at $(p,0)$, and let $f\colon X\to\mathbb P^1$
be the induced map. We have $f^*(0)=D+E$, where
$D=\operatorname{Bl}_p\mathbb T$ and $E\simeq\mathbb P^2$. Then $D$ is
very exceptional, and
$$
 D|_D=-e,\qquad e=D\cap E,\qquad \cO_E(D)=\cO_E(1).
$$
As $e$ is the only curve on $D$ and $D\cdot e=1$, the divisor $D>0$ has
nonnegative degree on all fibre curves, although $-e$ has negative degree
against a K\"ahler form on $D$. We also note that $f$ is
equidimensional, but a pullback from $\mathbb P^1$ would have equal
coefficients along $D$ and $E$.
\end{example}

\medskip\noindent\textbf{Good minimal models and Mori fibre spaces.}
We say a normal space is globally strongly $\Q$-factorial if every global
rank-one reflexive sheaf has an invertible reflexive power.

Let $\phi\colon X\dashrightarrow Y$ be a bimeromorphic map of normal
compact complex spaces, and let $B_Y=\phi_*B+\sum_E E$, where $E$ runs
over the prime divisors on $Y$ exceptional over $X$. We call $\phi$ a
\emph{good minimal model} of $(X,B)$ if $(Y,B_Y)$ is dlt, $Y$ is
globally strongly $\Q$-factorial, $K_Y+B_Y$ is semiample, and on a
common resolution $p\colon W\to X$, $q\colon W\to Y$ we have
\begin{equation}\label{comparison}
 p^*(K_X+B)=q^*(K_Y+B_Y)+F,\qquad F\geq0,\qquad q_*F=0,
\end{equation}
where $\Supp F$ contains the strict transform of every prime divisor on
$X$ contracted by $\phi$. By \eqref{comparison}, every divisor extracted
by $\phi$ has log discrepancy zero. A
bimeromorphic contraction is a bimeromorphic map which extracts no
divisors.

A \emph{Mori fibre space} of $(X,B)$ consists of a map $\phi$ as above
and a projective contraction $f\colon Y\to Z$ of an extremal ray, where
$Y$ and $Z$ are compact K\"ahler, $\dim Z<\dim Y$, and $-(K_Y+B_Y)$
is $f$-ample \cite[Definition~2.17]{HX}. We require $(Y,B_Y)$ to be
dlt, $Y$ to be globally strongly $\Q$-factorial, and \eqref{comparison}
with the same condition on $\Supp F$.

\medskip\noindent\textbf{Generalized pairs.}
For a generalized pair (g-pair) $(X/Z,B,\M)$ we follow
\cite[Definitions~2.7 and~2.9]{DHY}. On a log resolution $p\colon W\to X$
on which $\M$ descends, $M_W$ is nef over $Z$ and
$$
 c_1(K_W+B_W)+M_W=p^*(K_X+B+M_X),\qquad p_*B_W=B,
$$
where $B_W$ has snc support and $K_X+B+M_X\in H^{1,1}_{\mathrm{BC}}(X,\R)$.
The g-pair is klt (resp.\ lc) if $1-\coeff_E B_W>0$
(resp.\ $\geq0$) for each prime divisor $E$ on any such resolution.
We will use freely that on a minimal model over $Z$ the adjoint class is
nef over $Z$, and that the difference of the two pullbacks to a common
resolution is effective and exceptional over the model.

\medskip\noindent\textbf{The divisorial Zariski decomposition.}
Recall that a pseudoeffective real $(1,1)$-class $\alpha$ on a compact
complex manifold $X$ has the divisorial Zariski decomposition
\cite[\S3]{Bou}
$$
 \alpha=P+c_1(N_\sigma(\alpha)),
$$
where $P$ is modified nef, that is, $\nu_D(P)=0$ for every $D$.
The sum defining $N_\sigma(\alpha)$ is finite.
For divisor classes on a projective manifold, $\nu_D$
equals Nakayama's $\sigma_D$ \cite[Theorem~5.4 and
Appendix~A.1]{Bou}. We write $N_\sigma(L)=N_\sigma(c_1(L))$ for an
$\R$-line bundle $L$. On a normal compact complex space we define
$N_\sigma(\alpha)=p_*N_\sigma(p^*\alpha)$ for any resolution $p$
\cite[Definition~A.6 and Lemma~A.7]{DHY}.

\begin{lemma}\label{negative-part}
The negative part has the following properties.
\begin{enumerate}
\item Let $\beta$ be a pseudoeffective class on a compact complex
manifold and let $N=N_\sigma(\beta)$. If $E\geq0$ is supported on $N$,
then $N_\sigma(\beta+c_1(E))=N+E$. If $0\leq E\leq N$, then
$N_\sigma(\beta-c_1(E))=N-E$.
\end{enumerate}
Let $p\colon W\to X$ be a modification from a compact complex manifold
$W$ onto a normal compact complex space $X$, and let
$\alpha\in H^{1,1}_{\mathrm{BC}}(X,\R)$ be pseudoeffective.
\begin{enumerate}
\setcounter{enumi}{1}
\item For every effective $p$-exceptional $\R$-divisor $E$, we have
\begin{equation}\label{exceptional}
 N_\sigma(p^*\alpha+c_1(E))=N_\sigma(p^*\alpha)+E .
\end{equation}
\item If $N=N_\sigma(\alpha)$ is $\R$-Cartier and $\alpha-c_1(N)$ is
nef, then $N_\sigma(p^*\alpha)=p^*N$.
\end{enumerate}
\end{lemma}
\begin{proof}
For (1) see \cite[proof of Theorem~3.18]{Bou}, and for (2) see
\cite[Lemma~A.5]{DHY}. For (3) let $P=\alpha-c_1(N)$. Since $p^*P$ is
nef, $G=p^*N-N_\sigma(p^*\alpha)$ is effective (subadditivity).
We have $p_*G=0$. The positive part $p^*P+c_1(G)$ of $p^*\alpha$ has zero negative part.
By (2) that negative part is $G$, so $G=0$.
\end{proof}

We remark that in \eqref{comparison} we have $F=N_\sigma(p^*(K_X+B))$ by
\eqref{exceptional}, as $K_Y+B_Y$ is nef. In particular, the support of
$N_\sigma(K_X+B)$ is the union of the prime divisors contracted by the
good minimal model.

\begin{lemma}[Minimal models over the semiample model]\label{fixed-over-base}
Let $g\colon X\to Z$ be a surjection with connected fibres from a
normal compact complex space $X$ onto a normal projective variety $Z$.
Let $D$ be a pseudoeffective $\R$-Cartier divisor on $X$ and $H$ an
ample $\R$-divisor on $Z$ such that
$$
 D\sim_\R g^*H+E,\qquad E=N_\sigma(D).
$$
Then the following hold.
\begin{enumerate}
\item The vertical part of $E$ is very exceptional over $Z$. Suppose
moreover that $c_1(D)$ is nef on the fibres of $g$, and that $E$ is
vertical or $D|_F\sim_\R0$ for a very general fibre $F$ of $g$. Then
$E=0$ and $D\sim_\R g^*H$.
\item Suppose that $X$ is K\"ahler and globally strongly
$\Q$-factorial, and that $D=K_X+B$ for a dlt pair $(X,B)$. Then there
is a finite $(K_X+B)$-MMP over $Z$ ending with a good minimal model
$V$ of $(X,B)$ such that $K_V+B_V\sim_\R g_V^*H$, where
$g_V\colon V\to Z$ is the induced morphism.
\end{enumerate}
\end{lemma}
\begin{proof}
(1) We argue as in \cite[Lemmas~2.18 and 2.19]{Hu20}.
Let $Q\subset Z$ be a prime divisor. Let $R$ be the sum of the prime
divisors dominating $Q$, counted with their multiplicities in $g^*Q$
near a general point of $Q$. Suppose $\Supp R\subset\Supp E$.
Near a general point of $Q$ we have $g_*\cO_X(R)=\cO_Z(Q)$. For a very
ample $A$ and $m\gg0$, by global generation of
$g_*\cO_X(R)\otimes\cO_Z(mA)$ we can take an effective
$R'\sim R+mg^*A$ sharing no component with $R$.
For $0<\varepsilon\ll1$, the divisor $E-\varepsilon R+\varepsilon R'$ is
effective and $H-\varepsilon mA$ is ample. Hence
$N_\sigma(D)\le E-\varepsilon R+\varepsilon R'$, contrary to
$N_\sigma(D)=E$ along $R$.

If $D|_F\sim_\R0$, the effective divisor $E|_F$ has Gauduchon
degree zero on a resolution of $F$. Hence $E|_F=0$ and $E$ is vertical.
In either case $E$ is very exceptional and nef on the fibres, so
$E=0$ by Lemma~\ref{very-exceptional}(1).

(2) We argue as in the proof of \cite[Proposition~3.17]{Hu25}.
We take $kH\cdot C>2\dim X$ for every integral curve $C$ on $Z$ and run
a $(D+kg^*H)$-MMP with scaling of a sufficiently large K\"ahler class
$\theta$. By \cite[Theorem~5.16(1)]{HX}, for $0<t\ll1$ we reach a nef
model $V$ of $D+kg^*H+t\theta$. Indeed, for $S=\lfloor B\rfloor$ and
$0<\eta\ll t$, the pair $(X,B-\eta S)$ is klt and $t\theta+\eta c_1(S)$
is K\"ahler. By the length of extremal rays for $D+t\theta$
\cite[Theorem~1.3]{HX}, the MMP runs over $Z$, so each step is $D$-negative.

On a resolution, by Lemma~\ref{negative-part}(1) and (3) we have
$N_\sigma(D+kg^*H)=N_\sigma((1+k)D-kE)=E$.
By \cite[Proposition~3.5]{Bou}, for sufficiently small $t$ every
component of $E$ occurs in $N_\sigma(D+kg^*H+t\theta)$. By
\eqref{exceptional}, $V$ contracts $E$. We have $K_V+B_V\sim_\R g_V^*H$.
\end{proof}

\begin{lemma}\label{resolution-model}
Let $p\colon X'\to X$ be a log resolution of an lc pair $(X,B)$ and let
$\Delta=p_*^{-1}B+\Exc(p)_{\mathrm{red}}$. Suppose that a bimeromorphic
contraction $\psi\colon X'\dashrightarrow Y$ is a good minimal model of
$(X',\Delta)$. Then $\psi\circ p^{-1}$ is a good minimal model of
$(X,B)$. Moreover, $\psi$ contracts each $p$-exceptional prime divisor
of positive log discrepancy with respect to $(X,B)$.
\end{lemma}
\begin{proof}
In $K_{X'}+\Delta=p^*(K_X+B)+A$ the coefficients of $A$ are the
exceptional log discrepancies. On a common resolution $u\colon W\to X'$,
$v\colon W\to Y$, let $r=pu$ and let $G$ be as in \eqref{comparison} for
$\psi$. Since $\psi$ extracts no divisors, $B_Y=\psi_*\Delta$ and
$$
 r^*(K_X+B)-v^*(K_Y+B_Y)=G-u^*A .
$$
Since $u^*A-G$ is $r$-nef and $r_*(u^*A-G)\le0$, by
Lemma~\ref{very-exceptional}(1) we have $0\le G-u^*A\le G$. In
particular, $G-u^*A$ is $v$-exceptional. Its coefficients along the
strict transforms of prime divisors on $X$ are those of $G$, so
$\psi\circ p^{-1}$ is a good minimal model of $(X,B)$. Since
$u^*A\le G$, $\psi$ contracts every component of $A$.
\end{proof}
\section{Divisorial decompositions for boundaries with real coefficients}

On a smooth compact K\"ahler manifold $T$, let $\sum_{i=1}^h S_i$
be a reduced snc divisor with prime components. For $c\in[0,1]^h$
we write $\Delta(c)=\sum_i c_iS_i$.

\begin{lemma}\label{cube}
Suppose that $K_T+\Delta(c)$ is pseudoeffective for some
$c\in[0,1]^h$. Then the set of $c\in[0,1]^h$ such that
$K_T+\Delta(c)$ is pseudoeffective is a rational polytope. Moreover,
it admits a finite rational polyhedral decomposition such that, for each
prime divisor $D$ over $T$, the asymptotic order of $|K_T+\Delta(c)|_\Q$
along $D$ is affine on the rational points of each cell.
\end{lemma}
\begin{proof}
\emph{Step 1. Finite generation.} We realize the section spaces
as weight spaces of an adjoint ring on the fibre product
$\pi\colon\widetilde T\to T$ of the
$\mathbb P_T(\cO_T\oplus\cO_T(S_i))$, with tautological classes $H_i$.
On the $i$-th factor the quotients by $\cO_T(S_i)$,
$\cO_T\cdot(1,0)$ and $\cO_T\cdot(1,\sigma_i)$
($\sigma_i$ defines $S_i$) define sections of classes
$H_i-\pi^*S_i$, $H_i$ and $H_i$. The first section is disjoint from
the other two, whose equations are $y_i=0$ and $y_i+z_i=0$ near
$S_i=(z_i=0)$. Since the base divisor is snc, the sum of the pullbacks
is a reduced snc divisor $\Gamma$. From
$K_{\widetilde T/T}\sim\sum_i(-2H_i+\pi^*S_i)$ we get
$K_{\widetilde T}+\Gamma\sim\pi^*K_T+\sum_iH_i$. By pushforward to $T$,
$$
 H^0\big(\widetilde T,d(K_{\widetilde T}+\Gamma)\big)=
 \bigoplus_{0\leq\alpha_i\leq d}
 H^0\Big(T,dK_T+\sum_i\alpha_i S_i\Big)
 \qquad(d\geq0),
$$
where $\alpha_i$ is the degree in $\cO_T(S_i)$.

The adjoint divisor is pseudoeffective by
$K_{\widetilde T}+\Gamma\sim\pi^*(K_T+\sum_iS_i)+\sum_i(H_i-\pi^*S_i)$
(increase the boundary to the reduced divisor and add the first sections).
By \ref{gzd}, $\mu^*(K_{\widetilde T}+\Gamma)\sim_\Q P+N$ on a
modification $\mu$. For sufficiently divisible $r$, $rP$ is free and
$krN$ is the fixed part in degree $kr$ for every $k\ge1$ (minimal
multiplicity). Removing the fixed parts, we see that the $r$-th
Veronese subring is the finitely generated section ring of $rP$.

\emph{Step 2. Weights and orders of vanishing.} Let $w_j=(d_j,\alpha_j)$
be the weights of finitely many multihomogeneous generators of positive
degree for the Veronese subring. The fan of \cite[Lemma~2.1]{ELMNP}
on the cone generated by the $w_j$ cuts out a rational polyhedral
decomposition of the convex hull of the $\alpha_j/d_j$ in degree one.

For rational $c\in[0,1]^h$, membership in the hull is equivalent to
the existence of a monomial of weight $m(1,c)$ for some positive
integer $m$. Since monomials span the weight spaces, this means
$|K_T+\Delta(c)|_\Q\ne\varnothing$. By \ref{gzd}, this is equivalent
to pseudoeffectivity of $K_T+\Delta(c)$. Both sets are closed and
have dense rational points (approximate $c$ from above in $[0,1]$
in the pseudoeffective set), so they coincide.

For a prime divisor $D$ over $T$, let $e_j$ be the order of the $j$-th
generator along $D$. A monomial of weight $m(1,c)$ with exponents $n_j$
has normalized order $\sum_j(n_j/m)e_j$. For $c\in\Q^h$ in the hull,
a rational minimizer yields a monomial after clearing denominators.
Hence
\begin{equation}\label{orders}
 o_D(c)=\min\left\{\sum_j\lambda_j e_j\,
     \middle|\ \lambda_j\geq0,\quad\sum_j\lambda_j w_j=(1,c)\right\},
\end{equation}
By the choice of fan, the minimum is affine on each cell for any $D$.
\end{proof}
\begin{proposition}\label{realdecomp}
Let $(X,B)$ be an lc pair on a normal compact K\"ahler space. Suppose
that $K_X+B$ is pseudoeffective. Then there is a projective log
resolution $p\colon W\to X$ from a smooth compact K\"ahler manifold $W$
such that
$$
 p^*(K_X+B)\sim_\R P+N_\sigma(p^*(K_X+B))
$$
for a semiample $\R$-line bundle $P$. Moreover, if $B$ is a
$\Q$-divisor, then the decomposition holds up to $\Q$-linear
equivalence with $P$ a semiample $\Q$-line bundle.
\end{proposition}
\begin{proof}
\emph{Step 1. The log smooth case.} Suppose $(X,B)=(T,\Delta(c))$.
By Lemma~\ref{cube}, $c$ is in the relative interior of a cell with
vertices $v_j$. We write $c=\sum_jt_jv_j$ as a convex combination.
If $c\in\Q^h$, we take the $t_j$ rational. By Lemma~\ref{negative-part}(3),
on a projective log resolution $u\colon W\to T$ dominating the modifications
in \ref{gzd} at the vertices we have
$$
 u^*(K_T+\Delta(v_j))\sim_\Q P_j+N_j,\qquad
 N_j=N_\sigma(u^*(K_T+\Delta(v_j))),
$$
with $P_j$ semiample.

Let $D$ be a prime divisor on $W$ and $x$ a rational point of the cell.
The positive part has a free multiple on a common resolution of $W$
and the modification in \ref{gzd} for $x$. We read off
$\nu_D(u^*(K_T+\Delta(x)))=o_D(x)$ from the coefficient of the strict
transform of $D$ in the negative part. By Lemma~\ref{cube}, $o_D$ is affine
on rational points. By convexity in $x$, the minimal multiplicity
is continuous in the relative interior. We obtain at $c$
$$
 N_\sigma(u^*(K_T+\Delta(c)))=\sum_jt_jN_j,\qquad
 u^*(K_T+\Delta(c))-\sum_jt_jN_j\sim_\R\sum_jt_jP_j .
$$
The sum is semiample.

\emph{Step 2. The general case.} We apply Step~1 to $(T,\Delta)$, where
$f\colon T\to X$ is a projective log resolution and
$\Delta=f_*^{-1}B+\Exc(f)_{\rm red}$. Then $K_T+\Delta=f^*(K_X+B)+A$
with $A\ge0$ exceptional. By \eqref{exceptional}, the negative parts of
$u^*(K_T+\Delta)$ and $(fu)^*(K_X+B)$ differ by $u^*A$. We subtract
$u^*A$ to obtain the decomposition on $W$.
\end{proof}

\begin{proof}[Proof of Corollary~\ref{nef}]
Since $K_X+B$ is nef, by Proposition~\ref{realdecomp} and
Lemma~\ref{semiample} we have $p^*(K_X+B)\sim_\R h^*H$ on a resolution
$p\colon W\to X$. Here $h\colon W\to Z$ is a contraction and $H$ is
ample. Since the fibres of $p$ are connected and projective and their
curves have $h^*H$-degree zero, $h$ factors through a contraction
$g\colon X\to Z$. By injectivity of $p^*$ on Picard groups,
$K_X+B\sim_\R g^*H$. For the $\Q$-divisor case see Lemma~\ref{semiample}.
\end{proof}
\section{Termination with K\"ahler scaling}

We follow \cite[Lemma~3.14]{Hu25} to prove termination with scaling.

\begin{lemma}\label{interval}
Let $(X,B)$ be a dlt pair on a globally strongly $\Q$-factorial
compact K\"ahler space. Suppose that $K_X+B$ is pseudoeffective and
$\omega$ is a K\"ahler class such that $K_X+B+\omega$ is nef. Then there is a $\varepsilon>0$ such that in any $(K_X+B)$-MMP with
scaling of $\omega$, the model $X_i$ is a minimal model of $(X,B)$
whenever two consecutive scaling numbers satisfy
$0\leq\lambda_i<\lambda_{i-1}<\varepsilon$.
\end{lemma}
\begin{proof}
\emph{Step 1.} We construct a common minimal model $V$ for
$$
 p^*(K_X+B+s\omega)+(1-\delta)E,\qquad 0\le s\ll1 .
$$
Here $p\colon W\to X$ is the resolution of Proposition~\ref{realdecomp},
$E=K_W+p_*^{-1}B+\Exc(p)_{\rm red}-p^*(K_X+B)\ge0$, and $0<\delta\ll1$
is fixed. Let $\Delta=p_*^{-1}B+\Exc(p)_{\rm red}-\delta E$. By
\eqref{exceptional} and Lemma~\ref{semiample}, we have
$$
 K_W+\Delta=p^*(K_X+B)+(1-\delta)E\sim_\R g^*H+G,\qquad
 G=N_\sigma(K_W+\Delta),
$$
where $g\colon W\to Z$ is a contraction and $H$ is ample.

Let $S=\lfloor B\rfloor$ and take $\eta>0$ so small that
$\omega+\eta c_1(S)$ is K\"ahler. Since $(X,B-\eta S)$ is klt,
$\Supp p^*S$ contains the support of $\lfloor\Delta\rfloor$. For
$0<s_0\ll1$, the boundary $\Delta-s_0\eta p^*S$ and the nef and big
part $s_0p^*(\omega+\eta c_1(S))$ define a klt g-pair with adjoint
class $K_W+\Delta+s_0p^*\omega$. We fix $k\gg1$ as in Lemma~\ref{fixed-over-base}(2) and run an MMP
on $K_W+\Delta+s_0p^*\omega+kg^*H$ with K\"ahler scaling. For small $s_0$,
the MMP runs over $Z$ and ends with a model $V$ contracting $G$.
Thus $K_V+\Delta_V\sim_\R g_V^*H$. Denote the trace of $p^*\omega$ by
$\omega_V$. The nefness of $(1+k)g_V^*H+s_0\omega_V$ implies that
$K_V+\Delta_V+s\omega_V$ is nef for $0\le s\le s_0/(1+k)$.

On a common resolution $r\colon U\to W$, $w\colon U\to V$, let $F_s$
be the pullback difference for $K_W+\Delta+s p^*\omega$. We have
$F_0=r^*G$. By \eqref{exceptional}, $F_{s_0}$ is the negative part
of the pullback of $K_W+\Delta+s_0p^*\omega+kg^*H$. Since
$(1+k)g^*H+s_0p^*\omega$ is nef, $0\le F_{s_0}\le r^*G$.
By negativity the MMP is an isomorphism off $\Supp G$. On a common
log resolution that is an isomorphism there, $F_0$ has positive
coefficients along the $w$-exceptional divisors. Hence for
$0\le s\le s_0/(1+k)$ the affine difference $F_s$ is effective,
$w$-exceptional and positive along those divisors. Since $(W,\Delta)$
is log smooth and the nef part descends to $W$, this proves that $V$
is a minimal model throughout the interval.

\emph{Step 2.} Suppose $0\le\lambda_i<\lambda_{i-1}<s_0/(1+k)$.
For $s=\lambda_i,\lambda_{i-1}$, the adjoint class on $X_i$ is nef
by scaling, and its pullback difference from $W$ is effective and
exceptional over $X_i$ since $E$ is exceptional over $X_i$.
By comparing discrepancies on a common resolution
$v\colon U\to X_i$, $w\colon U\to V$ we get
$$
 v^*(K_{X_i}+B_i+s\omega_i)=w^*(K_V+\Delta_V+s\omega_V),
 \qquad s=\lambda_i,\lambda_{i-1}.
$$
By linearity the equality holds at $s=0$. So $K_{X_i}+B_i$ is nef
and $X_i$ is a minimal model.
\end{proof}

\begin{proof}[Proof of Theorem~\ref{scaling}]
We run a $(K_X+B)$-MMP with scaling of $\omega$. For
$0<\varepsilon<\lambda_i$, we take the klt g-pair as in the proof of
Lemma~\ref{fixed-over-base}(2), with adjoint class $K_X+B+\varepsilon\omega$.
Up to $X_i$ the steps form an MMP with scaling of $\omega$ for this
g-pair, with scaling number $\lambda_i-\varepsilon$ on $X_i$.
By \cite[Theorem~1.3 and
Proposition~4.3]{HX}, we can contract any ray for a $(K_{X_i}+B_i)$-step
at $\lambda_i$ and flip if necessary.

Suppose the MMP is infinite. If the $\lambda_i$ stay away from zero,
we can fix the klt g-pair and obtain termination by
\cite[Theorem~1.5 and proof of Theorem~5.16]{HX}. So $\lambda_i\to0$.
By nefness on $X_i$ and negativity, $K_X+B+\lambda_i\omega$ is
pseudoeffective for each $i$, as is its limit $K_X+B$. For some $i$
we have $0<\lambda_i<\lambda_{i-1}<\varepsilon$. By Lemma~\ref{interval},
$X_i$ is a minimal model, which leads to a contradiction.

If the final adjoint class is nef, it is semiample by
Corollary~\ref{nef}. By \eqref{comparison}, $K_X+B$ is pseudoeffective.
Otherwise the last step is a contraction $f\colon Y\to Z$ onto a
compact K\"ahler space with $\dim Z<\dim Y$ \cite[Proposition~4.3]{HX}.
On the contracted ray, $-(K_Y+B_Y)$ has the same degree as a positive
multiple of a K\"ahler class on $Y$. By \cite[Lemma~2.39]{HX}, $f$ is
projective and $-(K_Y+B_Y)$ is $f$-ample. Negativity along the MMP
gives \eqref{comparison} and the condition on $\Supp F$.
If $K_X+B$ were pseudoeffective, then
$K_X+B\sim_\R E\geq0$ by Proposition~\ref{realdecomp}. The transform
of $E$ would be both effective and antiample on a general fibre of $f$.
\end{proof}
\begin{proof}[Proof of Theorem~\ref{main}]
Suppose $K_X+B$ is not pseudoeffective. By
\cite[Theorems~2.27--2.28 and Corollary~2.29]{HX}, there is a projective
crepant modification $p\colon (X',B')\to(X,B)$ where $(X',B')$ is dlt
and $X'$ is globally strongly $\Q$-factorial. If $(X,B)$ is dlt, we
take $p$ small. Since $X'$ is compact K\"ahler and
$K_{X'}+B'=p^*(K_X+B)$ is not pseudoeffective, we obtain a Mori fibre
space $Y\to Z$ by Theorem~\ref{scaling}. By crepancy of $p$,
the pullback difference in \eqref{comparison} is the same
for $X\dashrightarrow Y$ and $X'\dashrightarrow Y$. A prime divisor on
$X$ contracted by $\phi$ has a strict transform on $X'$ contracted by
$X'\dashrightarrow Y$.

Suppose $K_X+B$ is pseudoeffective. Let $p\colon W\to X$ be a
projective log resolution of $(X,B)$. In the
dlt case we take the resolution with positive $p$-exceptional log
discrepancies.
Let $\Delta=p_*^{-1}B+\Exc(p)_{\rm red}$. Since $(X,B)$ is lc,
$K_W+\Delta=p^*(K_X+B)+A$ with $A\geq0$. Thus $K_W+\Delta$ is
pseudoeffective. Let $\omega$ be a K\"ahler class such that
$K_W+\Delta+\omega$ is nef. By Theorem~\ref{scaling}, a
$(K_W+\Delta)$-MMP with scaling of $\omega$ ends with a good minimal model.
By Lemma~\ref{resolution-model} we obtain a good minimal model of $(X,B)$.
In the dlt case the same lemma shows that every $p$-exceptional divisor
is contracted.
\end{proof}

\section*{Acknowledgments}

The first named author is supported by Overseas High-Level Young Talent Recruitment
Programs and is partially supported by the Chongqing Natural Science Foundation Innovation
and Development Joint Fund CSTB2023NSCQ-LZX0031. The second named author is partially
supported by the National Key Research and Development Program of China (Grant No.
2024YFA1013201), the Natural Science Foundation of Chongqing (Grant No.
CSTB2024NSCQ-LZX0040), and the Special Project of Chongqing Municipal Science and Technology
Bureau (Grant No.\ 2025CCZ015).

\end{document}